\documentclass{amsart}
\usepackage[PostScript=dvips]{diagrams}
\usepackage{amssymb}
\usepackage{amsmath}
\usepackage{amsthm}
\let\mathbb\mathbf 
\newcommand{\Z}{{\mathbb Z}}            
\newcommand{\Zhat}{\hat{\Z}}            
\newcommand{\Q}{{\mathbb Q}}            
\newcommand{\N}{{\mathbb N}}            
\newcommand{\G}{{\mathbb G}}            
\newcommand{\Zp}{{\mathbb Z_p}}         
\newcommand{\Zpstar}{{\mathbb Z^*_p}}           
\newcommand{\seteq}{\mathrel{:=}}
\newcommand{\Qp}{{\mathbb Q_p}}         
\newcommand{\Fp}{{\mathbf F}_p}         
\newcommand{\Kbar}{{\overline{K}}}              
\newcommand{\Qbar}{{\overline{\Q}}}
\newcommand{\m}{{\mathfrak{m}}}         

\newcommand{\Section}{\S}
\newcommand{\into}{\hookrightarrow}     
\renewcommand{\H}{{\mathcal H}}         

\newcommand{\T}{{\mathbf T}}            
\newcommand{\I}{{\mathfrak{I}}}         

\renewcommand{\P}{{\mathbf P}}          

\newcommand{\la}{\langle}               
\newcommand{\ra}{\rangle}               

\newcommand{\s}{\sigma}                 

\def\hlinefill{\leaders\hrule height 3pt depth -2.5pt\hfill}
\def\emrule{\thinspace\hbox to 0.75em{\hlinefill}\thinspace}

\DeclareMathOperator{\Gal}{{\rm Gal}}           
\newcommand{\GQ}{\Gal (\overline{\Q} / \Q )}    
\DeclareMathOperator{\GL}{{\mathbf{GL}}}                
\DeclareMathOperator{\WM}{{\rm WM}}             
\DeclareMathOperator{\WG}{{\rm WG}}             

\newcommand{\F}{{\mathbf F}}                    

\newtheorem{theorem}{Theorem}[section]
\newtheorem{conj}[theorem]{Conjecture}
\newtheorem{cor}[theorem]{Corollary}
\newtheorem{lemma}[theorem]{Lemma}
\newtheorem{prop}[theorem]{Proposition}
\newtheorem{example}[theorem]{Example}

\def\pf{{\sc Proof. }}
\def\qed{\hfill $\blacksquare$\smallskip}

\theoremstyle{remark}
\newtheorem{remark}[theorem]{Remark}
\newtheorem{defn}{Definition}

\newcounter{listcounter1}
\newenvironment{list1}{
  \begin{list}{\arabic{listcounter1}.\hfill}{
    \usecounter{listcounter1}
    \setlength{\leftmargin}{18pt}
    \setlength{\labelwidth}{18pt}
    \setlength{\labelsep}{0pt}
    \setlength{\topsep}{0pt}
  }
}{
  \end{list}
}

\begin{document}

\title[Galois Theory and Torsion Points on Curves]
        {Galois Theory and Torsion Points on Curves}


\author[M.~Baker]{Matthew H.~Baker}
\address{Department of Mathematics, Harvard University,
                Cambridge, MA 02138, USA}
\email{mbaker@math.harvard.edu}

\author[K.~Ribet]{Kenneth A.~Ribet}
\address{Department of Mathematics, University of California,
                Berkeley, CA 94720-3840, USA}
\email{ribet@math.berkeley.edu}

\thanks{The authors' research was partially supported by an NSF 
Postdoctoral Fellowship and by NSF grant
DMS-9970593.  The authors thank
Bjorn Poonen and Ron Fertig for providing
useful comments on earlier drafts of this manuscript.}


\maketitle

\section{Introduction}
\label{Introduction}

This paper surveys
Galois-theoretic techniques for studying torsion 
points on curves that have been developed in recent years 
by A.~Tamagawa and the present authors.  

\medskip

We begin with a brief history of the problem of determining the set of
points of a curve that map to torsion points of the curve's Jacobian.

\medskip

Let $K$ be a number field, and suppose that 
$X/K$ is an algebraic curve\footnote
{
By an algebraic curve, we mean a complete, 
nonsingular, and absolutely irreducible variety of dimension one 
over a field.
}
of genus $g \geq 2$.  Assume, furthermore, that $X$ is embedded 
in its Jacobian variety $J$ via
a $K$-rational Albanese map $i$;
thus
there is a 
$K$-rational divisor $D$ of degree~one on $X$
such that $i = i_D : X \into J$ is defined on $\Kbar$-valued points
by the rule $i(P) = [(P)-D]$, where $[\,\cdot\,]$ denotes 
the linear equivalence class of a divisor on $X$.  When $D$ is a 
$K$-rational point $P_0$, we often refer to $P_0$ as the {\em base point}
of the embedding $i_Q$.

\medskip

Let $T \seteq J(\Kbar)^{\rm tors}$ denote the torsion subgroup of $J(\Kbar)$.

\begin{theorem}
\label{MMconjecture}
The set $X(\Kbar) \cap T$ is finite. 
\end{theorem}

Theorem~\ref{MMconjecture} 
was stated as the
{\em Manin--Mumford conjecture\/} by S.~Lang in 1965.
In his article~\cite{Lang2}, Lang reduced this conjecture to a second
conjectural statement, which concerns
the action of Galois groups on
torsion points of abelian varieties over finitely generated fields.
This latter statement is still unproven, despite recent partial
progress by Serre,
Wintenberger (see \cite{Wintenberger}) and other authors.
The first proof of the Manin--Mumford conjecture was provided
by
M.~Raynaud
\cite{Raynaud1}, who combined
Galois-theoretic results on torsion points of $J$ with a subtle analysis
of the reductions mod $p^2$ of $X$ and $J$ for a suitable prime $p$.  
A second proof was given by R.~Coleman\footnote
{
The results of Tamagawa that we present in 
section \ref{tamagawa-rtpc} of this paper are closely related to Coleman's
work in \cite{RTPC}, although the methods are different.
}
in \cite{RTPC} using $p$-adic integration to
analyze the set of primes that
may ramify in the field generated by a torsion
point on $X$.  

\medskip

Raynaud also proved the following generalized
version of the Manin--Mumford conjecture (see \cite{Raynaud2}):

\begin{theorem}
\label{raynaudtheorem}
Let $K$ be a field of characteristic zero, and let $A/K$ be an
abelian variety.  Let $V$ be a subvariety of $A$ which is not the translate
by a torsion point of a positive-dimensional abelian subvariety $B$ of $A$.
Let $T \seteq A(\Kbar)^{\rm tors}$.  Then the set $T \cap V$ is not 
Zariski-dense in $V$.
\end{theorem}

One can use Theorem \ref{raynaudtheorem} to establish uniform bounds
for the cardinality of $X(\Kbar) \cap T$ as one varies the
Albanese embedding; see \cite{BakerPoonen} for details.

\medskip

It is also possible to generalize 
Theorem~\ref{raynaudtheorem} in several different
directions, replacing $T$ by the division group of any 
finitely generated subgroup of $A(\Kbar)$, or by any sequence of points in
$A(\Kbar)$ whose canonical height tends to zero.  
See \cite{Hindry}, \cite{McQuillan}, and \cite{Poonen1} for precise statements and further results.  

\medskip

In this paper, however, we focus on the original problem:
What can we say about the intersection $X(\Kbar) \cap T$ when $X$ is
a curve?  We are particularly interested in
explicit
determination of this intersection for particular classes of
curves.  We mention the following three results:

\medskip

1. ({\em Curves of genus 2}) B.~Poonen's paper \cite{Poonen2}
gives an algorithm (which has been implemented on a computer)
for determining the intersection $X(\Qbar)\cap T$
when $X/ \Q$ is a genus 2 curve embedded in its Jacobian 
using a Weierstrass point.  Poonen's method relies crucially on ideas 
of Buium \cite{Buium} and Coleman \cite{RTPC}.

\medskip

2. ({\em Fermat curves\/})
Suppose $X$ is the plane curve given by the equation $x^m + y^m = z^m$ for 
$m\geq 4$.  
The {\em cusps\/} of $X$ are the points 
$(x,y,z)\in X(\Qbar)$ such that $xyz = 0$.  
Rohrlich \cite{Rohrlich}
proved that the difference of 
two cusps is always torsion as an element of~$J$.  
Fix a cusp $c$ and embed $X$ in $J$ using $c$ as a base point. 
Coleman, Tamagawa, and Tzermias \cite{CTT} prove:
\begin{theorem}
\label{ctt-theorem}
The torsion points on $X$ in the embedding $i_c : X \into J$ are
precisely the cusps.
\end{theorem}

The proof of this theorem involves, among other things, 
Coleman's $p$-adic integration methods, 
complex multiplication theory, and results on class numbers of 
cyclotomic fields.

\medskip

3. ({\em Modular curves})
In \cite{Baker} and \cite{Tamagawa}, the authors independently prove a
conjecture of Coleman, Kaskel, and Ribet \cite{CKR}
concerning torsion points on 
the modular curve $X_0(p)$ in the cuspidal embedding.

Recall that a curve $X/K$ of genus $g\geq 2$ over a field $K$ is
{\em hyperelliptic\/} if 
there exists a degree 2 map $f : X \to \P^1$
defined over $\Kbar$.  Such a map, if it exists, is necessarily unique
(up to an automorphism of $\P^1$), and
the ramification points of $f$ are called the 
{\em hyperelliptic branch points}.

The Coleman--Kaskel--Ribet conjecture is the following statement.

\begin{theorem}
\label{ckr-theorem}
Let $p\geq 23$ be a prime number, and let $X$ be the modular curve
$X_0(p)$.
Let $H$ be the set of
hyperelliptic branch points on $X$ when $X$ is hyperelliptic and $p\neq 37$,
and otherwise let $H = \emptyset$.
Then the set of torsion points on $X$ in the embedding 
$i_\infty : X \into J$ is precisely $\{ 0, \infty \} \cup H$.
\end{theorem}

%

Note that the condition $p \geq 23$ in the statement of
Theorem~\ref{ckr-theorem} is equivalent to the genus of $X_0(p)$ being at least 2.

\medskip

We do not discuss results (1) or (2) further in this paper, but
we will say much more about the modular curves $X_0(p)$, and we
give a complete proof of Theorem \ref{ckr-theorem} in section 
\ref{CKR-section}.  

\begin{remark}
It is easy to obtain results
similar to Theorem~\ref{ckr-theorem} for $X_0(mp)$ or $X_1(mp)$ with $p\geq 23$ prime and
$m$ arbitrary by utilizing the natural maps $X_1(mp) \to X_0(mp) \to X_0(p)$.  
See \cite[Proposition~4.1]{Baker} for details.
\end{remark}

\begin{remark}
Though the proof of Theorem~\ref{ckr-theorem} we give in
  this paper is simpler than the previously published ones, 
  it still relies upon a number of deep results, e.g. Grothendieck's
  semistable reduction theorem, Mazur's detailed study of the
  arithmetic of $X_0(p)$ and $J_0(p)$, and the second author's level-lowering theorem.
\end{remark}

\medskip

Here is a brief outline of the contents of this paper.  In section 2 we 
discuss what it means for an element of a module to be 
``almost fixed'' by a group action, and we prove some elementary lemmas
about such elements.  We then show how these 
ideas can be combined with a result of Serre to give a simple proof of
the Manin--Mumford conjecture.  In section 3, we study torsion 
points on Abelian varieties which are almost fixed by the action
of an inertia group.  This is done, following Tamagawa, 
in the abstract setting of ``ordinary semistable'' and 
``ordinary good'' modules.  In section 4, the abstract algebraic 
manipulations of section 3 are placed in a geometric context, with
Theorems \ref{TMT} and \ref{OS-ramification}
as the reward.  In section 5, we discuss the proof of Theorem~
\ref{ckr-theorem}.  We attempt to give references for all of the 
facts we use about modular curves and their Jacobians.  The material in
section 5 relies on section 2 up through and including 
Lemma 2.7, and on section 3 up through
Theorem 3.6, so the reader who is only interested in reading the proof of
Theorem~\ref{ckr-theorem} can skip section 4 and the other parts of
sections 2 and 3.  In order to preserve the flow of the paper, a few
results quoted in the body of the paper are relegated to appendices.

%
%

\medskip
\medskip

{\small
{\bf Acknowledgements:} 
%

We include fairly detailed proofs of all results presented in this paper
in order to keep the exposition reasonably self-contained.  
However, a number of the
proofs in this paper can also be found in \cite{Baker}, \cite{Tamagawa},
and \cite{Kim-Ribet}.  
All results in sections
\ref{OS-modules} and \ref{tamagawa-rtpc}, except for Proposition
\ref{splitting}, are due to Tamagawa, and appear in his 
paper \cite{Tamagawa}.  
However, most of the proofs in section \ref{OS-modules} are new.
The proof of 
Theorem \ref{CKR-conj} which we give combines elements from both
\cite{Baker} and \cite{Tamagawa}.  



Commutative diagrams in this paper were designed using Paul Taylor's
Commutative Diagrams in \TeX\ package.
}

\section{Almost Rational Points and the Manin--Mumford Conjecture}
\label{AlmostRational}


In this section $K$ is a field and $X/K$ is an algebraic curve 
of genus at least 2.

\medskip

The results of this section and the next are motivated by the following
simple observation, which
plays a 
key role in the proof of the Coleman--Kaskel--Ribet conjecture.  

\begin{lemma}
\label{keylemma}
Suppose $X$ is embedded in its Jacobian $J$ via a 
$K$-rational Albanese map $i_D$.
Let $P\in X(\Kbar)$; if $X$ is hyperelliptic, assume that $P$ is 
not a hyperelliptic branch point.
Suppose that there exist $g,h \in \Gal(\Kbar / K)$ such that $gP + hP = 2P$
in $J$.  Then $gP = hP = P$.
\end{lemma}

\pf To be pedantic, we write $Q = i_D (P)$, so that $P$ is a point on
$X$ and $Q = [(P)-(D)]$ is its image in the Jacobian of $X$.
We are given that $gQ + hQ = 2Q$ in $J$, so that the degree-zero divisors
$(gP) - (gD) + (hP) - (hD)$ and $2(P) - 2(D)$ are linearly equivalent.
Since $D$ is $K$-rational,
it follows that the divisors $(gP) + (hP)$ and $2(P)$ on $X$ are
linearly equivalent, so that
there exists a rational function $f$ on $X$ whose
divisor is $(gP) + (hP) - 2(P)$.  Since $P$ is not a hyperelliptic branch
point, $f$ must be constant, so that $gP = hP = P$, as desired.\footnote
{
Notice how we are using
that $J$ is both the Albanese and Picard
variety for $X$.  The interplay between the two properties of $J$ 
lies behind many of the geometric results discussed in this paper.
}
\qed

\medskip

Lemma \ref{keylemma} suggests the following definition.

\medskip

\begin{defn}
Let $G$ be a group, and let $M$ be a $\Z [G]$-module.  
An element $P$ of $M$ is
{\em almost fixed\/} (by~$G$)
if $(g+h-2)P=0$ with $g,h \in G$ implies 
that $(g-1)P = (h-1)P = 0$.

The module $M$ is
almost fixed if $(g+h-2)M = 0$ with
$g,h \in G$ implies that $(g-1)M = (h-1)M = 0$.
\end{defn}

\begin{remark}
If $G=G_K$ is the absolute Galois group of a field $K$, 
we will often use the term {\em almost rational\/} instead of
almost fixed.
\end{remark}

We will be particularly interested in the set of almost rational 
{\em torsion points\/} of $M$.

\begin{example}
\label{x+y=2}
The set of almost rational torsion points of $\G_m(\Qbar)$ 
is~$\mu_6$, the group of sixth roots of unity.  
\end{example}

The proof is left as an exercise for the reader
(or see \cite[Lemma 3.14]{Baker}).

\medskip

We now prove some elementary lemmas 
concerning almost fixed elements and almost fixed modules.

\begin{lemma}
\label{basiclemma}
Let $P$ be an almost fixed element of the $\Z [G]$-module $M$.
\begin{list1}
\item If $\s \in G$, then $\s P$ is almost fixed.
\item If $g \in G$ and $(g - 1)^2 P = 0$, then $(g-1)P = 0$.
\end{list1}
\end{lemma}

\pf 
For the first part, notice that if $(g+h-2)\s P = 0$, then 
\[
(\s^{-1} g \s + \s^{-1} h \s - 2)P = 0, 
\]
which implies that
$(\s^{-1} g \s - 1)P = (\s^{-1} h \s - 1)P = 0$.  Therefore both $g$ and
$h$ fix $\s P$, as desired.

\medskip

For the second statement, 
we are given that $(g^2 - 2g + 1)P = 0$.  Multiplying on the left by
$g^{-1}$, we find that $(g + g^{-1} - 2)P = 0$, and therefore $(g-1)P=0$
by the definition of ``almost fixed.''
\qed

\begin{lemma}
\label{generators}
Let $M$ be a $\Z [G]$-module.  If $M$ is generated by almost fixed
elements, then $M$ is almost fixed.
\end{lemma}

\pf Let $P_1,\ldots,P_k$ be almost fixed elements that generate $M$ as a 
$\Z [G]$-module, and let $g,h$ be elements of $G$ such that 
$(g+h-2)M=0$.  Then $(g+h-2)(\s P_i) = 0$ for all $\sigma\in G$ and all
$i = 1,\ldots,k$.  By Lemma \ref{basiclemma}, each $\s P_i$ is almost
fixed, and therefore both $g$ and $h$ fix all of the $\s P_i$.  As
the $\s P_i$ generate $M$ as a $\Z$-module, it follows that both 
$g$ and $h$ fix every element of $M$.  Therefore $M$ is almost fixed.
\qed

\begin{remark}
It is not true that if $M$ is almost fixed then every element of $M$ is 
almost fixed.  For example, let $M$ be the 2-dimensional
$(\Z/5\Z)$-vector space $(\Z/5\Z)^2$, and let $G = (\Z/4\Z)$ act on
$M$ by sending a generator to 
$A := \left( \begin{array}{ll} 0 & 1 \\ -1 & 0 \\ \end{array} \right).$
A short computation shows that $M$ is almost fixed,
but the vector $v := \left[ \begin{array}{l} 2 \\ 1 \\ \end{array}
\right]$ 
is not, since $Av + A^2v = 2v$ but $Av \neq v$.
\end{remark}

Let us return to the geometric situation of Lemma
\ref{keylemma}, 
so that $K$ is a field, $G_K = \Gal (K^{\rm sep} / K)$ 
is the absolute Galois group of $K$,
and $X/K$ is a curve of genus $g\geq 2$,
embedded in its Jacobian $J$ via a $K$-rational Albanese map.  


\medskip


If $P\in X(\Kbar)$, then following A.~Tamagawa,
we say that the pair $(X,P)$ is {\em exceptional\/} if $X$ is hyperelliptic 
and $P$ is a 
hyperelliptic branch point on $X$.

\medskip

The following is a reformulation of Lemma \ref{keylemma} using our new
terminology.

\begin{lemma}
\label{curvelemma}
Let $P$ be a $\Kbar$-valued point of $X$.
Then either $(X,P)$ is exceptional, or $P$ is almost rational.
\end{lemma}


We illustrate the usefulness of the notion of almost rationality by
presenting a short proof
of the Manin--Mumford conjecture.

\medskip

The proof exploits the following deep result\footnote
{
This result was presented in Serre's Coll{\`e}ge de France lectures 
(1985--1986),
but the proof has not yet been published.  The main theorems of
\cite{Hindry} and \cite{McQuillan} both depend on this result.
}
due to Serre.


\begin{theorem}
\label{Serre-theorem}
Let $K$ be a finitely generated field of characteristic zero.  Let 
$A/K$ be an abelian variety of dimension $g$, and let $\rho : G_K \to 
\GL_{2g}(\Zhat)$ denote the Galois representation arising from the adelic
Tate module of $A$.  
Let $\Zhat^* \subset \GL_{2g}(\Zhat)$ denote the subgroup of homotheties.
Then the group $\Zhat^* / \left(\rho(G_K) \cap \Zhat^*\right)$ 
has finite exponent.
\end{theorem}

We will also need the following lemma (compare with Example \ref{x+y=2}):

\begin{lemma}
\label{countingpoints}
Let $e$ be a positive integer.  Then there is a positive constant $C(e)$
such that for all integers $m > C(e)$, 
there exist $x,y \in (\Z / m\Z)^*$ such that
$x^e,y^e \neq 1$ but $x^e + y^e = 2$.
\end{lemma}

\pf
By the Chinese remainder theorem, it suffices to consider the case
where $m = p^k$ is a prime power.  

\medskip

If $k = 1$, we want to look at $\Fp$-rational points on the projective
curve $C$ defined by $x^e + y^e = 2$.
By the Weil bounds, $\# C(\Fp) = p+1 + O(\sqrt{p})$.  Since
the number of points $(x,y) \in C(\Fp)$ with one of $x^e,y^e$ being
0 or 1 is at most $(e+1)^2$, the result follows in this case.  

\medskip

Finally, suppose $k\geq 2$.  If $p>e$,
Hensel's lemma guarantees the existence of $x,y \in \Z / p^k \Z$
such that $x^e = 1 + p^{k-1}, \; y^e = 1 - p^{k-1}$.  Since $x^e y^e = 1$,
we have $x,y \in (\Z / m\Z)^*$.
\qed

We can now prove the following finiteness result:

\begin{theorem}
\label{Ribet-theorem}
Let $K$ be a finitely generated field of characteristic zero, 
and let $A/K$ be an abelian variety.  
Then the set of almost rational torsion points on $A$ is finite\footnote
{
See also \cite{Calegari}, in which the author
classifies almost rational torsion points on
semistable elliptic curves over $\Q$.
}.
\end{theorem}

\pf
By Theorem \ref{Serre-theorem}, there exists a positive integer $e$ 
such that the group $\Zhat^* /\left( \rho(G_K)\cap \Zhat^*\right)$
has exponent $e$.  
Let $P$ be a torsion point on $A$ of order $m > C(e)$.  By Lemma
\ref{countingpoints}, 
there exist $x,y \in (\Z / m\Z)^*$ such that
$x^e,y^e \neq 1$ but $x^e + y^e = 2$.  
Since $(\Zhat^*)^e \subseteq \rho(G_K)\cap \Zhat^*$,
we can choose $g,h \in \GQ$ such that
$g,h$ act on $A[m]$ as $x^e$ and $y^e$, respectively.  Then
$(g+h-2)P = 0$ but neither $g$ nor $h$ fixes $P$, so $P$ is not
almost rational.  It follows that the set of almost rational torsion points
on $A$ is finite.
\qed

The Manin--Mumford conjecture follows easily from \ref{Ribet-theorem}:

\begin{cor}
Let $K$ be as above, and let $X$ be a curve of genus at least $2$, 
embedded in its Jacobian $J$ by an Albanese map.  
Then the set of torsion points on $X$ is finite.
\end{cor}

\pf The set of hyperelliptic branch points on $X$ is finite, as is
the set of almost rational torsion points on $J$.  The result therefore
follows from Lemma \ref{curvelemma}.
\qed

%

\section{Ordinary semistable and almost unramified modules}
\label{OS-modules}

In this section, $R$ will denote the ring of integers in a finite 
{\em unramified\/} extension $K$ of $\Qp$, where $p$ is an {\em odd\/} prime.%
\footnote{For the case $p=2$, see Tamagawa's paper \cite{Tamagawa}.}



We will denote by $I$ the inertia subgroup of $G \seteq \Gal(\Kbar / K)$, and
by $I^{\rm tame}$ the inertia subgroup of $\Gal(K^{\rm tame}/K)$, where
$K^{\rm tame}$ is the maximal tamely ramified extension of $K$.  Recall 
that the group $I^{\rm wild} \seteq \Gal(\Kbar / K^{\rm tame})$ is a pro-$p$
group, and that $I^{\rm tame}$ is canonically isomorphic to 
the group $\varprojlim \F_{p^n}^*$, 
where the transition maps are given by taking norms.  

For each $n\geq 1$, we denote by $I(n)$ the (normal) subgroup 
of $I$ fixing all of the $p^n$th roots of unity in $\Kbar$, and we let
$I(\infty)$ be the intersection of $I(n)$ for all natural numbers $n$, so
that $I(\infty)$ is the subgroup of $I$ fixing all $p$-power roots of unity.

\medskip

The motivation for the results in this section comes from the 
following observation:

\begin{lemma}
\label{trivial1}
Let $X/K$ be a curve of genus at least 2, embedded in its Jacobian $J$
via a $K$-rational Albanese map.  Suppose that $J$ is semistable, that 
$P\in X(\Kbar)$ is a 
torsion point of order prime to $p$, and that $(X,P)$ is not exceptional.
Then $I$ fixes $P$, i.e., $P$ is unramified.
\end{lemma}

\pf
Grothendieck showed in \cite[Proposition 3.5]{Groth} 
that if $A/K$ is a semistable abelian variety and $P\in A(\Kbar)$ 
has order prime to $p$, then $(\sigma - 1)^2 P = 0$ for all $\sigma\in I$.
Therefore in our situation we have 
\[
\s P + \s^{-1}P - 2P = \s^{-1}(\sigma -1)^2 P = 0
\]
for all $\s\in I$.  The result now follows from Lemma \ref{curvelemma}.
\qed


We now make some definitions.

\medskip

\begin{defn}
Let $M$ be a $\Z [I]$-module.
An element $P\in M$ (resp.\ $M$ itself) is
{\em almost unramified\/} if $P$ (resp.\ $M$) is almost fixed with
respect to the action of $I$.
\end{defn}

In other words, $M$ is almost unramified
if and only if whenever $(g+h-2)M = 0$ with $g,h \in I$, we have
$(g-1)M = (h-1)M = 0$. 

\medskip


\begin{defn}
A finite $\Z [I]$-module
$M$ is {\em ordinary semistable\/} if there exists
an exact sequence of $\Z [I]$-modules
\begin{equation}
\label{filtration}
0 \to M' \to M \to M'' \to 0
\end{equation}
such that:

\medskip

(i) $I$ acts on $M'$ via\footnote
{
If $N$ is a torsion abelian group, then $N$ is naturally a $\Zhat$-module.  
Also, the inertia group $I$ comes equipped with a cyclotomic
character $\chi: I \to \Zhat^*$.  It therefore makes sense 
to say that a torsion $I$-module $N$ is cyclotomic: this means that 
$\sigma n = \chi(\sigma)n$ for all $\s\in I$ and $n\in N$.
Note that if $N$ is cyclotomic and has order prime to $p$, 
then $I$ acts trivially on $N$, and that in general if $N$ is 
cyclotomic, then $I$ will act on $N$ through its abelian 
quotient $I/I(\infty)$.
}
the cyclotomic character $\chi$

(ii) $I$ acts trivially on $M''$.

\end{defn}
\medskip
For each finite $\Z [I]$-module $M$, there is a unique decomposition 
$M = M_p \oplus M_{\text{non-}p}$, where
$M_p$ has $p$-power order and $M_{\text{non-}p}$ has order prime to $p$.
Using this notation, we have the following definition.

\begin{defn}
A finite $\Z [I]$-module $M$
is {\em ordinary good\/} if it is ordinary semistable and,
in addition, $I$ acts trivially on $M_{\text{non-}p}$.  
\end{defn}


The definitions of ordinary good and ordinary semistable modules 
are motivated by the following:

\begin{defn}
An abelian variety
$A/K$
has {\em ordinary
semistable reduction\/} if
the connected component of the closed fiber of the N{\'e}ron model of $A$ 
over $R$ is an extension of an ordinary abelian variety by a torus.
\end{defn}

\begin{theorem}
\label{sga7}
Let $A$ be an abelian variety over $K$ and let $n$ be a positive
integer.
\begin{list1}
\item If $A$ has good ordinary reduction over $R$, then $A[n]$ is an
ordinary good $\Z [I]$-module.
\item If $A$ has ordinary semistable reduction over $R$,
then $A[n]$ is an ordinary semistable 
$\Z [I]$-module.
\end{list1}
\end{theorem}

\pf This is a consequence of Grothendieck's study in SGA7 of 
Galois actions on torsion points of semistable abelian varieties.
See \cite{Tamagawa} for details and precise references.
\qed

As a prototype of results to come, we have the following lemma (compare
with Lemma \ref{trivial1}):

\begin{lemma}
\label{triviallemma}
Suppose $M$ is a finite ordinary semistable and almost unramified 
$\Z [I]$-module of order prime to
$p$.  Then $I$ acts trivially on $M$.
\end{lemma}

\pf
Let $M'$ and $M''$ be as in the definition of ``ordinary semistable''. 
Since the order of $M'$ is prime to $p$ and $M'$ is cyclotomic, $I$ acts
trivially on $M'$.  Since $I$ acts trivially on $M''$ as well, it
follows that $(g-1)^2 M = 0$ for all $g \in I$.
But since $M$ is also almost unramified, Lemma \ref{basiclemma} tells us
that $(g - 1) M = 0$, so that $I$ acts trivially on $M$, as desired.
\qed

With an eye toward applying the results of this section 
to the study of ramified torsion 
points on curves, we now undertake an investigation of 
modules that are both almost unramified and ordinary semistable.

\begin{lemma}
\label{abelian}
If $M$ is a finite ordinary semistable and almost unramified 
$\Z [I]$-module, then $I(\infty)$ acts
trivially on $M$.  Therefore, the action of
$I$ on $M$ factors through its abelian quotient 
$I/I(\infty) \cong \Zpstar$.  
\end{lemma}

\pf 
Since $I(\infty)$ acts trivially on both $M'$ and $M''$ in 
the filtration (\ref{filtration}) coming from the definition of 
``ordinary semistable,'' it follows that
$(g-1)^2 M=0$ for all $g\in I(\infty)$.
That $M$ is almost unramified
then implies, by Lemma \ref{basiclemma}, 
that $I(\infty)$ acts trivially on $M$.
\qed

\begin{prop}
\label{chiprop}
Let $M$ be a finite ordinary semistable and almost unramified 
$\Z [I]$-module.  Let $p^m$ be the
order of $M_p$, and let $g,h$ be elements of $I$ such that 
$\chi(g) + \chi(h) \equiv 2$ mod $p^m$.  
Then $(g+h-2)M_p = 0$.
\end{prop}


\pf 
Let $0 \to M' \to M \to M'' \to 0$ be the filtration of $M$ given by
(\ref{filtration}), let $M'_p = M' \cap M_p$, and let $M''_p$ be the
image of $M_p$ in $M''$ under the given surjection.  

Then we have an exact sequence
\begin{equation}
\label{p-filtration}
0 \to M'_p \to M_p \to M''_p \to 0
\end{equation}
of modules of $p$-power order which again satisfies properties (i) and (ii)
in the definition of ``ordinary semistable.''

Since $I$ acts on $M'_p$ via the cyclotomic character, and since we are
assuming that $p>2$, the subgroup $(M'_p)^I$ of inertia invariants in 
$M'_p$ must be zero.

\medskip

The identity $\chi(g) + \chi(h) = 2$ mod $p^m$ implies that
$\alpha \seteq g + h - 2$ kills both $M_p'$ and $M_p''$.  
Therefore $\alpha$ acts on 
$M_p$ via a homomorphism $\phi: M_p'' \to M_p'$.
Since the action of $I$ on $M$ is abelian by Lemma \ref{abelian}, $\phi$
is a homomorphism of $I$-modules, and therefore (since $I$ acts trivially
on $M''$) the image of $\phi$ is contained in $(M_p')^I = 0$.
It follows that $\alpha$ kills $M_p$, as desired.
\qed


\begin{theorem}
\label{tametheorem}
Let $M$ be a finite ordinary semistable and almost unramified 
$\Z [I]$-module.  

\medskip

(1) The group $I(1)$ acts trivially on $M$.  

(2) If $p \geq 5$ and $M$ is ordinary good, then $I$ acts trivially on $M$.
\end{theorem}

\pf Let $g\in I(1)$.
Since $\chi(g)$ is 1 mod $p$ and $\chi : I \to \Zpstar$ is surjective, 
we can find $h$ in $I(1)$ such that 
$\chi(g) + \chi(h) = 2$ in $\Zp$.
By Proposition \ref{chiprop}, 
$(g+h-2)M_p =0$, where $M_p$ again
denotes the $p$-primary part
of $M$.  Also, by Lemma \ref{abelian} we know that $I(\infty)$ acts
trivially on $M$, from which it follows by Lemma \ref{pro-p} and the
definition of ``ordinary semistable'' that the action of the pro-$p$ group 
$I(1)/I(\infty)$ on $M_{\text{non-}p}$ is trivial.
Therefore $(g+h-2)M = 0$. As $M$ is almost unramified, 
it follows that $(g-1)M = 0$.  This proves part (1) of the theorem.

\medskip

To prove part (2), assume that $p\geq 5$ and that $M$ is ordinary good.
Then $g \in I$ acts
trivially on $M$ whenever we can solve the equation 
$\chi(g) + \chi(h) = 2$ in $\Zp$, i.e., whenever 
$\chi(g)$ is not 2 mod $p$.  Thus the set of $g \in I$ acting
trivially on $M$ forms a subgroup of $I$ whose image $H$ in $I/I(1)$ 
contains at least $p-2$ elements.
Since $I/I(1)$ has order $p-1$ and $p\geq 5$, it follows 
that $H = I/I(1)$.  Therefore $I$ acts trivially on $M$.  
\qed


We take a moment to remind the reader of our running assumption 
that $p$ is odd.

\begin{prop}
\label{splitting}
Let $M$ be a finite ordinary semistable and almost unramified 
$\Z [I]$-module.
Then $M'_p$ is killed by $p$ and $M_p = M'_p \oplus (M_p)^I$.
\end{prop}



\pf  By Theorem \ref{tametheorem}, $I(1)$ acts trivially on $M$.
Since $I$ acts via the $p$-adic cyclotomic character $\chi$ on $M'_p$, 
it follows that $pM'_p=0$.  

\medskip

To prove the second statement, consider the exact sequence 
\[
0 \to M'_p \to M_p \to M''_p \to 0 
\]
given by (\ref{p-filtration}).  
As we have already seen, $(M'_p)^I=0$.  
As $I$ acts on $M'_p$ through its abelian
quotient $I/I(1)$, we can apply Sah's lemma (Lemma \ref{Sah})
to an element $g$ of $I$ such that $\chi(g)$ is 2 mod $p$, and we see
that $H^1(I,M'_p)=0$.
Therefore the natural map $(M_p)^I \to (M''_p)^I = M''_p$ 
is an isomorphism, which is equivalent to the desired statement that
$M_p = M'_p \oplus (M_p)^I$.
\qed

\begin{cor}
\label{splittingcor}
In addition to the hypotheses of the proposition, suppose we are given 
an element $g\in I$ and an integer $r\in\Z$ such that 
$\chi(g) \equiv -r$ (mod $p$).  Then
$(g + r)(g + g^{-1} - 2) M = 0$.
\end{cor}

\pf By Proposition \ref{splitting}, we have 
\[
M = M'_p \oplus (M_p)^I \oplus M_{\text{non-}p}
\]
and $pM'_p = 0$, so it follows from the definition of 
``ordinary semistable'' that $$(g+r)(g-1)^2 M = 0.$$  Therefore
\[
g^{-1} (g+r)(g^2 - 2g + 1)M = 0,
\]
which gives the desired result.
\qed

\begin{cor}
\label{stabilizer}
Suppose $M$ is a finite cyclic ordinary semistable and almost
unramified $\Z [I]$-module, and that $P$ is a generator.
If $I$ acts nontrivially on $M$, then the 
group of elements of $I$ that fix $P$ is precisely $I(1)$.
\end{cor}

\pf We know by Theorem \ref{tametheorem} that $I(1)$ acts trivially on $M$.
Since $I$ acts nontrivially on $M$ but
trivially on $M_{\text{non-}p}$ (by Lemma \ref{triviallemma}), 
$p$ must divide the order of $M$.
We then see from Proposition \ref{splitting} 
that $p$ exactly divides the order of $M'_p$, and that
$I$ acts on $M'_p$ via the mod $p$ cyclotomic character.  
In particular, we can use Proposition \ref{splitting} to 
write $P$ as $x + y$, where $x \in M'_p$ and $y \in M^I$.
We must have $x \neq 0$, or else $I$ would fix $\s P$ for all $\s \in I$
and therefore act trivially on $M$.
It follows that
$(g - 1) P = (\chi(g)-1)x \neq 0$ for all $g\in I$ such that 
$\chi(g) \not\equiv 1$ (mod $p$), i.e., for all
$g \in I - I(1)$.
\qed

\section{Ramified torsion points on curves}
\label{tamagawa-rtpc}

As in the previous section, $K$ denotes a finite unramified extension of
$\Qp$, with $p\neq 2$.

Throughout this section, $X$ will denote a curve over $K$, 
embedded in its Jacobian $J$ via a $K$-rational Albanese map.

\medskip

In this section, we apply the results of section \ref{OS-modules} to
the study of torsion points on $X$.  The idea, due to Tamagawa, 
is to use elements of the
inertia group $I$ which act nontrivially on a torsion point $P\in X(\Kbar)$
to produce rational functions on $X$ of small degree.

\medskip

We first recall some basic facts about algebraic curves which can be
found, for example, in \cite[III.5]{FK}.

\medskip

If $P \in X(\Kbar)$, we denote by $\WM(P)$ the {\em Weierstrass monoid\/} at
$P$ consisting of all nonnegative integers $m$ such that there exists a 
rational function on $X$ of degree exactly $m$ having no poles outside $P$.
It is clear from the definition that $0 \in \WM(P)$, and that 
if $a,b \in \WM(P)$ then $a+b \in \WM(P)$, so that $\WM(P)$ is indeed a 
monoid.

\medskip
Let $\N$ denote the monoid $\{ 0,1,2,\ldots \}$ of nonnegative integers,
together with the operation of addition.
The complement of $\WM(P)$ in $\N$,
which we denote by $\WG(P)$, is called
the set of {\em Weierstrass gaps\/} at~$P$.   
It follows from the Riemann--Roch theorem
that $\WG(P)$ has exactly $g$ elements.  
A point $P$ on $X$ is called a {\em Weierstrass point\/} if 
there exists $m \in \WM(P)$ such that $1\leq m\leq g$, or equivalently, if 
$\WG(P) \neq \{ 1,2,\ldots,g \}$.  It is well
known that a curve of genus
$g\geq 2$ has at most $g^3 - g$ Weierstrass points.

\medskip

We now investigate the implications of the results of the previous section
for ramified torsion points on curves.

Part (2a) of the following theorem was originally proved
by Coleman using $p$-adic integration techniques.  The rest
of the theorem is due to Tamagawa.

\begin{theorem}
\label{TMT}
Let $X$ be a curve over $K$ whose Jacobian $J$ has ordinary
semistable reduction, and suppose $X$ is embedded in $J$ using a 
$K$-rational point.  

Let $P$ be a torsion point on $X$.  Then:

\smallskip

(1) The group $I(1)$ fixes $P$.

(2a) If $p \geq 5$ and $J$ has good ordinary reduction, then 
$P$ is unramified.

(2b) If $p=3$ and $J$ has good ordinary reduction, then either $P$ is
unramified or $3 \in \WM(P)$.
\end{theorem}

\pf 
When $(X,P)$ is exceptional, the result follows from 
Proposition \ref{exceptional}.  So we may assume that $(X,P)$ is not
exceptional.

\medskip

By Theorem \ref{sga7}, the $\Z [I]$-submodule $M$ of $J$ generated by $P$ is 
ordinary semistable, and is ordinary good 
when $J$ has good ordinary reduction.
Since $(X,P)$ is not exceptional, 
it follows from Lemmas \ref{curvelemma} and \ref{generators} that 
$M$ is almost unramified.  Parts (1) and (2a) therefore follow from 
Theorem \ref{tametheorem}.

For part (2b), note that if 
$\s \in I$ does not fix $P$, then $\s P - P$ has order $p$
in $J$ by Proposition \ref{splitting}.  Therefore the divisor $p(\s P) -
p(P)$ is principal.
\qed

\begin{prop}
\label{prop-nongap}
Suppose $J$ has ordinary semistable reduction, and
let $P$ be a torsion point of $J$ lying on $X$ which is ramified at $p$.  
Assume also that $(X,P)$ is not exceptional.
Let $r$ be a positive
integer
such that $r \not\equiv 0,1,$ or $-1$ 
(mod $p$).  Then 
the integer $2r-1$ lies in~$\WM(P)$; i.e., 
there exists a rational function of degree $2r-1$ on $X$ with no poles
outside $P$.
\end{prop}

\pf Let $M$ be the $\Z [I]$-module generated by $P$.  
Then as in the proof of Theorem \ref{TMT}, $M$ is ordinary semistable 
and almost unramified.
By hypothesis, $I$ acts nontrivially on $M$.
Also, by Corollary \ref{stabilizer}, $\sigma P \neq P$ for all $\sigma \in I$
such that $\chi(\sigma) \not\equiv 1$ (mod $p$).
Since $\chi$ is surjective, 
given any positive integer $r$ such that $r\not\equiv 0$ (mod $p$),
we can find $\sigma \in I$ such that 
$\chi(\sigma) \equiv -r$ (mod $p$).  
If in addition $r \not\equiv 1$ or $-1$ (mod $p$), 
then $\sigma^2 P \neq P$.  By
Corollary \ref{splittingcor}, we also know that
$(\sigma + \sigma^{-1} - 2)(\sigma + r) P = 0$ in $J$.
Multiplying this expression out, we find that
there exists a rational function $f$ on $X$ whose divisor is
\[
(\sigma^2 P) + (r-2)(\sigma P) + r(\sigma^{-1} P) - (2r-1)(P).
\]
The proposition now follows from the fact that the 
degree of $f$ is $2r-1$, since $P$ does not equal $\sigma^{-1} P$, 
$\sigma P$, or $\sigma^2 P$.
\qed

\medskip

The following is one of the main theorems of
Tamagawa~\cite{Tamagawa}.

\begin{theorem}
\label{OS-ramification}
Assume that $J$ has ordinary semistable reduction, 
that $(X,P)$ is not exceptional, and that 
$P$ is a ramified torsion point on $X$.  Then:
\begin{list1}  
\item If $p\geq 5$, then $g \leq 4$.  
\item If $p \geq 7$, then $g \leq 3$.  
\item If $p \geq 29$, then $g \leq 2$.
\end{list1}
\end{theorem}

\pf Suppose, for example, that $p \geq 5$.  
Taking $r = 2,3$ in Proposition \ref{prop-nongap}, we see that
$3,5 \in \WM(P)$.  By Lemma \ref{postage}, it follows that
$\WG(P) \subseteq \{ 1,2,4,7 \}$, and therefore $g\leq 4$.
Similarly, if $p\geq 7$, then taking $r=4$ we find that $7$ is also in
$\WM(P)$, and therefore $\WG(P) \subseteq \{ 1,2,4 \}$, so that $g \leq 3$.
Finally, suppose $p\geq 29$ and $g=3$.  We know from
Corollary \ref{stabilizer}
that the stabilizer of $P$ in $I$ is precisely $I(1)$.  Therefore
the set $\{ \sigma P \; | \; \sigma \in I \}$ has $p-1 \geq 28$ elements.
Since $3 \in \WG(P)$, $P$ must be a Weierstrass point, and therefore 
all of the points $\sigma P$ with $\sigma \in I$ must be Weierstrass 
points.
Since there are at most $g^3 - g = 24$ Weierstrass points on $X$, this
is a contradiction.
\qed

%

\medskip

We conclude this section with an intriguing open problem.
The following conjecture was made by R.~Coleman~\cite{RTPC}:

\medskip

\begin{conj}
Let $p\geq 5$ be a prime number, and suppose that $K/\Qp$ 
is an unramified finite extension.
Let $X/K$ be a curve of genus $g\geq 2$, embedded in its Jacobian
via a $K$-rational Albanese map.
Suppose furthermore that $X$ has good reduction over $K$.
Then every torsion point $P\in X(\Kbar)$ is unramified.
\end{conj}

In \cite{RTPC}, Coleman proved this
conjecture in the following
cases:

(i) $X$ has ordinary reduction

(ii) $X$ has superspecial reduction

(iii) $p > 2g$.

\medskip

The hypotheses of the conjecture are necessary\emrule see 
\cite[Appendix]{Baker-thesis} for an example.

On the other hand, Theorem \ref{OS-ramification} shows that with a few more
restrictions on the prime $p$, the conclusion of the 
conjecture remains true if $X$ merely has ordinary 
semistable reduction over $K$.  It would be interesting to try to use the
Galois-theoretic methods surveyed in this paper to prove additional
cases of Coleman's conjecture.

\section{Torsion points on modular curves}
\label{CKR-section}

In this section, we use the results of section \ref{OS-modules} to
give a short proof of the Coleman--Kaskel--Ribet conjecture.

\medskip


We first recall some
facts about the modular curves $X_0(p)$, for
which a basic reference is Mazur~\cite{Mazur} (see also~\cite{Mazur2}).

\medskip

Fix a prime number $p\geq 5$.
The modular curve $X_0(p)$ is a
compactified coarse moduli space for degree-$p$ isogenies
between elliptic curves.

As a Riemann surface, $X_0(p)$ can be
thought of as the quotient of the complex upper half plane $\H$ by the
action of the group $\Gamma_0(p)$,
suitably compactified by adding the two cusps $0$ and $\infty$.
As an algebraic curve, $X_0(p)$ is defined over $\Q$ and the cusps $0$
and $\infty$ are $\Q$-rational points.

\medskip

From now on we assume that $p\geq 23$, which is equivalent to assuming
that the genus $g$ of $X_0(p)$ is at least 2.

To simplify notation, we let $X \seteq X_0(p)$ and 
$J \seteq J_0(p)$.

\medskip


There is an involution $w_p$ of $X$, called the Atkin--Lehner
involution, which interchanges $0$ and $\infty$.
We note that $w_p$ always has fixed points (\cite[\Section2]{Ogg}).

The quotient of $X$ by $w_p$ will be denoted by $X_0^+(p)$, or simply
$X^+$.  Its genus will be denoted by $g^+$.  

For $p \geq 23$, we have $g^+ = 0$ if and only if $p \in 
\{ 23, 29, 31, 41, 47, 53, 71 \}$.

It is known (see \cite{Ogg}) that $X$ is hyperelliptic if and only if either
$g^+ = 0$ or $p=37$.

\medskip

For each $Q \in X(\Qbar)$, we can define an embedding $i_Q$
of $X$ into $J$
by sending $P\in X(\Qbar)$ to the linear equivalence class 
of the degree-zero divisor $[(P) - (Q)]$.

We call $i_\infty$ the {\em standard embedding\/} of $X$ into $J$, and we 
let $T_\infty$ be the set of torsion points on $X$ in the standard embedding.

\medskip

We now recall the Coleman--Kaskel--Ribet conjecture (see
Theorem~\ref{ckr-theorem}).

\begin{theorem}
\label{CKR-conj}
For all prime numbers $p\geq 23$,
\[
T_\infty = \left\{ \begin{array}
{l@{\quad{\rm if}\quad}l}
\{ 0,\infty \} & g^+>0 \\
\{ 0,\infty \} \cup \{ {\rm hyperelliptic \; branch \; points} \} & g^+=0.
\end{array} \right.
\]
\end{theorem}

Before we can prove the conjecture, we need to review some
more facts about $X$ and $J$.  We begin with some definitions and
elementary facts, all of which can be found in \cite{Mazur}.

\medskip

The {\em cuspidal subgroup\/} $C$ of $J$ is
the cyclic subgroup
of $J$ generated by the class of the degree-zero 
divisor $(0)-(\infty)$ on $X$.  

The {\em Shimura subgroup\/} $\Sigma$ of $J$
is the kernel of the map $J_0(p) \to J_1(p)$ induced via
Picard functoriality from the natural map $X_1(p) \to X_0(p)$.

Both $C$ and $\Sigma$ have order $n \seteq (p-1)/({\gcd}(p-1,12))$.

\medskip

The endomorphism ring of $J_{\Qbar}$ 
contains (and in fact equals) the
{\em Hecke algebra\/} $\T$ generated by $w_p$ and by the Hecke operators 
$T_l$, with $l$ prime and different from~$p$.

The {\em Eisenstein ideal\/} is the ideal $\I$ of~$\T$
generated by $w_p + 1$ and
the differences
$T_l - (l+1)$ for $l\neq p$.  
A maximal ideal $\m$ of $\T$ is {\em Eisenstein\/} if it contains~$\I$.

\medskip
The subgroup 
\[
J[\I] \seteq 
\{\, P \in J(\Qbar) \, | \, tP = 0 \hbox{ for all }
t \in \I \,\}
\]
contains both $C$ and $\Sigma$.  We list below some additional properties
of this subgroup which we will need\emrule see \cite{Csirik} for a 
complete picture of $J[\I]$ as a Galois module.


\medskip

In addition to the above definitions and relatively simple facts, 
the proof of Theorem \ref{CKR-conj} will also require 
the following ten more difficult
facts about $X$ and $J$.  For the reader's benefit, we provide 
references and/or sketch the proofs for each of these facts.

\newcounter{mazurfacts}
\begin{list}{\bf 
\arabic{mazurfacts}:}{\usecounter{mazurfacts}}%
\def\jtem#1\par{\item{\sl #1}\hfill\break}
\jtem 
$J$ has good reduction outside $p$, and has purely toric (hence
ordinary semistable) reduction at $p$.
 
This is due to Igusa and Deligne--Rapoport.  See 
\cite[Theorem A.1]{Mazur} for a discussion and references.

\jtem
$J(\Q)^{\rm tors} = C$.

This is \cite[Theorem 1]{Mazur}.

\jtem
If $P \in X(\Q) \cap J(\Q)^{\rm tors}$, then $P \in \{ 0,\infty \}$.

When $p\neq 37,43,67,163$, this is a consequence of the fact
that, by \cite[Theorem 7.1]{Mazur}, $X(\Q) = \{ 0,\infty \}$.  
For the four exceptional cases,
see \cite[Proof of Proposition 1.2]{CKR}.

\jtem
The natural map $\Z\to\T/\I$ induces an isomorphism
$\Z / n\Z\approx\T/\I$.

This is \cite[II, Proposition 9.7]{Mazur}.
\jtem
$J[\I]$ is a free $\T / \I$-module of rank 2.
%

This follows from the analysis in \cite[Ch.~II, \Section16--18]{Mazur}, as
noted in \cite[\Section 3]{Ribet1}.
\jtem
The set of torsion points of $J(\Qbar)$ that are
unramified at all primes
above $p$ is precisely $J[\I]$.

This is \cite[Proposition 3.3]{Ribet1}. 
\jtem
Let $M$ be a finite torsion $\T [\GQ]$-submodule of $J(\Qbar)$, and let  
$V$ be a Jordan--H{\"o}lder 
factor of $M$.  Let $\m$ be the maximal ideal in $\T$ that annihilates 
$V$
and
consider $V$ as a representation of $\GQ$ over
the field~$\T / \m$. 
Then
if $\m$ is Eisenstein, 
then $V$ is one-dimensional and isomorphic to either $\Z / l\Z$ or $\mu_l$,
where $l$ is the characteristic
of~$\T / \m$.
If $\m$ is not Eisenstein, 
then $V$ is isomorphic to the standard
two-dimensional irreducible representation $\rho_{\m} : \GQ \to \GL_2(k)$
attached to $\m$.

See \cite[Chapter II]{Mazur} for a proof, and \cite[Theorem 2.1]{Ribet1}
for a discussion of the proof.
\jtem
Suppose $l \mid n$, and let $I$ be an inertia subgroup at $l$ of $\GQ$. 
If $M$ is a $\Z [I]$-module such that 
$M \subseteq J[\I]$, then $M$ is ordinary good.
%

This follows from Fact 7, together with results of Oort and Tate
on finite flat group schemes of prime order.  See \cite[Proposition 2.3,
$(v)\Rightarrow (i)$]{Tamagawa} for details.
\jtem
If $\m \mid p$, then $\rho_{\m}$ is not {\em finite\/} at $p$ in the
sense of
\cite[\Section 2.8]{Serre2}.

This is a consequence of Mazur's level-lowering theorem (see
\cite[Theorem 1.1]{Ribet2}), since
if $\rho_{\m}$ were finite at $p$, it would have to be modular of level 1,
which is impossible.  
\jtem
If $\m \mid p$, then $\rho_{\m}(I)$ is non-abelian for every inertia
group $I$ of $\GQ$ at $p$.

We sketch an argument similar to the one given in 
\cite[\Section 4, (1-2)]{Tamagawa}:
Let $M$ be the $\T / \m [I]$-module giving rise to $\rho_{\m}$.  Then
$M$ is ordinary semistable as a $\Z [I]$-module, so that
$M$ has a filtration $0 \to M' \to M \to M'' \to 0$ in which
$I$ acts trivially on $M''$ and on $M'$ via $\chi$.  As in
the proof of Proposition \ref{splitting},
if the action of $I$ on $M$ is abelian, then Sah's lemma (Lemma \ref{Sah})
shows that
$M = M' \oplus M''$, and therefore $M$ is finite at $p$.  This contradicts 
Fact~9.
\end{list}
\medskip
\noindent
{\em Proof of Theorem \ref{CKR-conj}:}

Let $P$ be a point of $X$ such that $i_\infty (P)$ is torsion.

When $(X,P)$ is exceptional, the result follows from 
\cite[Proposition 1.1]{CKR}\footnote
{
We briefly recall the argument.  
For $p\neq 37$, the fact that the hyperelliptic branch points are torsion 
points in the embedding $i_\infty$ follows directly from 
the fact that in those cases, $w_p$ coincides with 
the hyperelliptic involution.  For if $P$ is fixed by $w_p$, then since
the hyperelliptic involution acts as $-1$ on $J$, we have
\[
2[(P) - (\infty)] = [(P) - (\infty)] + [w_p(P) - w_p(0)] = [(0)-(\infty)],
\]
which is torsion.  The case $p=37$ is more complicated, and 
follows from explicit calculations found in \cite[\Section 5]{MaS}.
}.  
So we will assume from now on that $(X,P)$ is not exceptional.

\medskip

By Fact 3,
it is enough to prove that $P$ is defined over $\Q$.

\medskip

{\em Claim 1:} $P$ is unramified at $p$.

\medskip

\pf

Let $I$ be an inertia subgroup at $p$ of $\GQ$.
Since $J$ has ordinary semistable reduction at $p$ by Fact 1, 
and since $(X,P)$ is not exceptional, it follows from 
Theorem \ref{tametheorem}
that $I(1)$ fixes $P$.  Applying the same argument to every
conjugate of $P$, we see that $I$ acts on 
the $\T[\GQ]$-module $M$ generated by $P$ through its abelian quotient
$I/I(1)$.

\medskip

If $p$ divides the order of $M$, then
$I$ acts through an abelian quotient on some 
Jordan--H{\"o}lder factor $V$ of $M$ 
associated to a maximal ideal $\m$ of residue characteristic $p$.
But Fact 10 tells us that the action of $I$ on $V$ is 
necessarily non-abelian, a contradiction.

\medskip

Therefore $M$ has order prime to $p$.  By Lemma \ref{triviallemma}, 
it follows that $I$ acts trivially on $M$.  Since this is true for all
inertia groups $I$ at $p$, it follows that $P$ is unramified at $p$.
\qed

\medskip

{\em Claim 2:} $i_\infty(P) \in J[\I]$.  

\medskip

\pf
This follows from Fact 6 and Claim 1.
\qed

\medskip

{\em Claim 3:} If $P$ is not a cusp then $g^+ = 0$.

\medskip

\pf Let $Q \seteq i_\infty (P)$.
Since $w_p$ interchanges the two cusps on $X$, there is a 
unique cusp on $X^+$, which we also call $\infty$.
So the fiber of the degree two map $\pi: X\to X^+$ over 
$\infty$ is just $\{0,\infty \}$.
Let $J^+$ be the Picard (Jacobian) variety of $X^+$.  
The fact that $J^+$ is also the Albanese variety of $X^+$ implies 
there is a commutative diagram
\begin{diagram}
X               & \rTo^{i_{\infty}}     & J             \\
\dTo^\pi        &               & \dTo_{\pi_*}   \\
X^+             & \rTo^{i_{\infty}}     & J^+            \\
\end{diagram}

If $\pi^* : J^+ \to J$ denotes the map induced by Picard functoriality,
then the composite map $\pi^* \circ \pi_* : J\to J$ is the map $1+w_p$.
Also, $\pi^*$ is injective; this is a consequence 
(see \cite[Lemma 6]{BakerPoonen}) of the 
fact that $w_p$ has fixed points.
Since $\I$ contains $1+w_p$, it follows that if 
$Q\in J[\I]$, then $Q$ is sent to zero
under the projection $\pi_*$.

\medskip

Therefore, when $g^+ > 0$ (so that the map $i_\infty : X^+ \to J^+$ is
an embedding), we have $P=0$ or $P=\infty$ as desired.  
\qed

\medskip

{\em Claim 4:} $P$ is unramified at $2$ and $3$.

\medskip

\pf By Claim 3, we may assume that
$g^+ = 0$, i.e., that $p$ belongs to the set of prime numbers 
$\{ \,23, 29, 31, 41, 47, 53, 71\, \}$.
An explicit calculation shows that $3 \nmid n$, and that
$2 \mid n$ if and only if $p=41$.

So by Claim 2 and Fact 7, we are reduced to the case $p=41$, 
where we have $n=10$.  We need to show in this case that $P$ is unramified
at $2$.
Since $4 \nmid n$, it follows from Fact 5 that
$M_2$ is killed by 2. 

\medskip

Let $I$ be an inertia group of $\GQ$ at 2, and suppose 
there exists $\s \in I$ such that $\s P \neq P$.  
Since $J$ has good reduction at $2$, $I$ acts trivially on 
$M_{\text{non-}2}$, so
$\s Q - Q \in M_2$, and therefore 
$2(\s Q - Q) = 0$.  It follows that the divisor
$2(\s P) - 2(P)$ is principal on $X$, so $(X,P)$ is exceptional, a 
contradiction.
\qed

\medskip

{\em Claim 5:} $P$ is defined over $\Q$.

\medskip
\pf
By Fact 7(i) and Claim 2, $P$ is
unramified at all primes $l \nmid n$.
It suffices to show that $P$ is unramified at all $l\geq 5$ 
such that $l \mid n$.
Fix such a prime $l$ and an inertia group $I$ at $l$ in $\GQ$.
Let $M$ be the $\Z [I]$-submodule of $J[\I]$ 
generated by $Q$.  
By Fact 8, $M$ is an ordinary good $\Z [I]$-module.
%
Also, since $(X,P)$ is not exceptional, it follows from Lemma~
\ref{curvelemma} that $M$ is almost unramified.
Theorem \ref{tametheorem} then implies that $I$ acts trivially on $M$,
as desired.
\qed

This concludes the proof of Theorem \ref{CKR-conj}.

\medskip

For generalizations to torsion points on $X$ in noncuspidal Albanese 
embeddings into $J$,
and to certain other modular curves, plus an application to 
Mordell--Weil ranks, see \cite[\Section 4]{Baker}.

\medskip
\medskip

\appendix

\section{Some elementary algebraic results}

For the sake of completeness, we give the statements and 
proofs of some elementary algebraic results used in this paper.

\begin{lemma}
\label{pro-p}
Let $G$ be a group, and let $M$ be a
finite $\Z [G]$-module of order prime to $p$.  Suppose that the action
of $G$ on $M$ factors through a finite $p$-group $G'$, and that
$(g-1)^2 = 0$ for all $g \in G$.  Then $G$ acts trivially
on $M$.
\end{lemma}

\pf Let $q = p^k$ be the order of $G'$, and let $g \in G$.  Then
\[
0 = (g^q - 1)M = ([1 + (g-1)]^q - 1)M = q(g-1)M
\]
by the binomial theorem.  Since $M$ has order prime to $p$, it follows that
$(g-1)M=0$.
\qed

The following elementary result from group cohomology is
known as Sah's lemma.  Our proof is adapted from \cite[Lemma 8.8.1]{LangDG}.

\begin{lemma}[Sah's lemma]
\label{Sah}
Let $G$ be a group, let $M$ be a $G$-module, and let $g$ be in the
center of $G$.  Then $H^1(G,M)$ is killed by the endomorphism $x \mapsto 
gx - x$ of $M$.  In particular, if this endomorphism is an automorphism, 
then $H^1(G,M) = 0$.
\end{lemma}

\pf Let $f:G\to M$ be a 1-cocycle.  Then for all $h\in G$,
\[
f(h) = f(ghg^{-1}) = f(g) + gf(hg^{-1}) = f(g) + g[f(h) + hf(g^{-1})].
\]

Therefore 
\[
(g-1)f(h) = gf(h) - f(h) = -f(g) - ghf(g^{-1}) = -f(g) - hgf(g^{-1}).
\]
But the cocycle condition implies that $f(1)=0$, so
\[
0 = f(1) = f(gg^{-1}) = f(g) + gf(g^{-1})
\]
and therefore $(g-1)f(h) = (h-1)f(g)$, so that $(g-1)f$ is a coboundary.
\qed

\medskip

Recall that a {\em monoid\/} is a a set $S$ together with an associative
composition law on $S$ and an identity element $e \in S$.

We denote by $\N$ the monoid consisting of all nonnegative integers.

If $a_1,\ldots,a_k \in \N$, we denote by $\la a_1,\ldots,a_k \ra$ 
the monoid 
\[
\{ n_1 a_1 + \cdots + n_k a_k \; | \; n_i \in \N \}.
\]
It is the smallest submonoid of $\N$ containing $a_1,\ldots,a_k$.

\medskip


The following result is sometimes called the ``postage stamp lemma'':

\begin{lemma}
\label{postage}
If $a,b$ are relatively prime positive integers and $m$ is any
integer such that $m\geq (a-1)(b-1)$, then $m\in \la a,b \ra$.
\end{lemma}

\pf Since no two of the $b$ integers $m-ar$ ($0 \leq r \leq b-1$) are
congruent modulo $b$, one of them must be divisible by $b$, say
$m - ar_0 = bs_0$.  As
\[
bs_0 = m - ar_0 \geq (a-1)(b-1) - a(b-1) = -b + 1,
\]
we must have $s_0 \geq 0$, so that $m \in \la a,b \ra$ as claimed.
\qed

\section{The exceptional case}


In this appendix,
$K$ denotes a finite unramified extension of
$\Qp$ with $p\neq 2$, and $X/K$ is a curve of genus at least 2, 
embedded in its Jacobian $J$ via a $K$-rational Albanese map.

\medskip

The following result, which is essentially \cite[Proposition 3.1]{Tamagawa},
was used in the proof of Theorem~\ref{TMT}. 


\begin{prop}
\label{exceptional}
Suppose $J$ has ordinary semistable reduction.
Let $P \in X(\Kbar)$ be a torsion point, and suppose $(X,P)$ is exceptional.
Then:

\smallskip

(1a) $\s^2 P = P$ for all $\s \in I$.  

(1b) The group $I(1)$ fixes $P$.

(2) If $J$ has good ordinary reduction, then $P$ is unramified.

\end{prop}


\pf
Let $M$ be the $\Z [I]$-submodule of $J$ generated by $P$.  Since $P$ is a 
Weierstrass point on $X$, so is $\s P$, and therefore the divisors
$2(P)$ and $2(\s P)$ on $X$ are linearly equivalent for all $\s \in I$.
It follows that $2(\s - 1)P = 0$ in $M$.  Applying the same argument to
every conjugate of $P$, we see that $I$ acts trivially on $2M$.  In 
particular, since $p$ is odd, $(\s - 1) M_p = 0$ for all $\s \in I$.  

\medskip

Also note that by Lemma \ref{pro-p}, $I^{\rm wild}$ acts trivially on 
$M_{\text{non-}p}$, and therefore
$I$ acts on $M$ through its quotient $I^{\rm tame}$.

\medskip

If $J$ has ordinary good reduction, then
$(\s - 1) M_{\text{non-}p} = 0$ for all $\s \in I$ and therefore $I$
acts trivially on $M$ as desired.

\medskip

In general, since $M$ is ordinary semistable, we have 
$(\s - 1)^2 M_{\text{non-}p} = 0$ for all $\s \in I$.  
Since $(\s - 1)M_p = 0$ as
well, we see that in fact $(\s - 1)^2 M = 0$ for all $\s \in I$.  Adding this
to the relation $2(\s - 1)M=0$, we find that $(\s^2 - 1)M=0$ for all
$\s \in I$.  This proves (1a).  Statement (1b) now follows from the fact
that $I(1)$ is contained in the subgroup of $I$ topologically generated
by $\{ \s^2 \; | \; \s \in I \}$.  
Explicitly: $I$ acts on $M$
through a finite quotient $I'$ of $I^{\rm tame}$ isomorphic to $\F_{p^n}^*$
for some $n\geq 1$.  
The image of $\s$ in $I'$ has
norm 1 in $\Fp$ if and only if $\s\in I(1)$.
The result now follows from the fact that an element of
$\F_{p^n}^*$ is a square if and only if its norm to $\F_{p}^*$ is a square.
\qed


\end{document}